\documentclass[11pt,reqno]{amsart}
\usepackage{amsthm}
\usepackage{amssymb}
\usepackage{graphicx}
\usepackage{latexsym}
\usepackage{multicol}
\usepackage{verbatim,enumerate}
\usepackage{accents}
\usepackage[usenames]{color}
\usepackage[colorlinks=true, linkcolor=blue, citecolor=blue, urlcolor=blue]{hyperref}
\usepackage{hyperref}
\usepackage{amsmath, amscd}

\theoremstyle{definition}
\newtheorem*{cor}{Corollary}%[section]
\newtheorem*{lem}{Lemma}
\newtheorem*{prop}{Proposition}

\newtheorem*{thm}{Theorem}

\theoremstyle{definition}
\newtheorem*{defn}{Definition}
\newtheorem*{conj}{Conjecture}
\newtheorem*{example}{Example}

\newtheorem*{rem}{Remark}

\newtheorem*{theom}{Theorem}

\newcommand{\nc}{\newcommand}

\newcounter{cnt}
\numberwithin{equation}{section}

\makeatletter
\def\section{\def\@secnumfont{\mdseries}\@startsection{section}{1}%
  \z@{.7\linespacing\@plus\linespacing}{.5\linespacing}%
  {\normalfont\scshape\centering}}
\def\subsection{\def\@secnumfont{\bfseries}\@startsection{subsection}{2}%
  {\parindent}{.5\linespacing\@plus.7\linespacing}{-.5em}%
  {\normalfont\bfseries}}
\makeatother

 \nc{\Hom}{\operatorname{Hom}}
  \nc{\mode}{\operatorname{mod}}
\nc{\End}{\operatorname{End}} \nc{\wh}[1]{\widehat{#1}} \nc{\Ext}{\operatorname{Ext}} \nc{\ch}{\text{ch}} \nc{\ev}{\operatorname{ev}}
 \makeatletter
\def\mydggeometry{\makeatletter\dg@YGRID=1\dg@XGRID=20\unitlength=0.003pt\makeatother}
\makeatother \theoremstyle{remark}

\numberwithin{equation}{section}
\let\bwdg\bigwedge
\def\bigwedge{{\textstyle\bwdg}}

\newcommand{\wt}{\operatorname{wt}}

\newcommand{\rnc}{\renewcommand}

\nc{\cal}{\mathcal} \nc{\goth}{\mathfrak} \rnc{\bold}{\mathbf}

\nc\bomega{{\mbox{\boldmath $\omega$}}} \nc\bpsi{{\mbox{\boldmath $\Psi$}}}
 \nc\balpha{{\mbox{\boldmath $\alpha$}}}
 \nc\bpi{{\mbox{\boldmath $\pi$}}}
 \nc\bvpi{{\mbox{\boldmath $\varpi$}}}

\newcommand{\lie}[1]{\mathfrak{#1}}

\newcommand\Lg{\mathfrak{g}}
\newcommand\Lh{\mathfrak{h}}

\newcommand\Ln{\mathfrak{n}}

\nc{\C}{\mathbb C }
\nc{\Z}{\mathbb Z }
\nc{\N}{\mathbb N }
\nc{\R}{\mathbb R }
\nc{\Q}{\mathbb Q }
\nc\blambda{{\mbox{\boldmath $\lambda$}}}
\nc\bmu{{\mbox{\boldmath $\mu$}}}
\nc\bsigma{{\mbox{\boldmath $\sigma$}}}

\newcommand{\sh}{\operatorname{sh}}
\newcommand{\PBW}{\operatorname{PBW}}

\newcommand{\gr}{\operatorname{gr}}

\DeclareMathOperator{\D}{D}
\DeclareMathOperator{\W}{W}
\DeclareMathOperator{\V}{V}
\DeclareMathOperator{\T}{T}
\DeclareMathOperator{\loc}{loc}
\begin{document}

\title[Fusion products and toroidal algebras]{Fusion products and toroidal algebras}

\author{Deniz Kus}
\address{Deniz Kus:\newline
Mathematisches Institut, Universit\"at zu K\"oln, Germany}
\email{dkus@math.uni-koeln.de}

\author{Peter Littelmann}
\address{Peter Littelmann:\newline
Mathematisches Institut, Universit\"at zu K\"oln, Germany}
\email{peter.littelmann@math.uni-koeln.de}

\thanks{This work has been partially supported by the DFG Priority Program SPP 1388 \textit{Representation Theory} and the “SFB/TR 12--\textit{Symmetries and
Universality in Mesoscopic Systems}”.}

\subjclass[2010]{}
\begin{abstract}
We study the category of finite--dimensional bi--graded representations of toroidal current algebras associated to finite--dimensional complex simple Lie algebras. Using the theory of graded representations for current algebras, we construct in different ways objects in that category and prove them to be isomorphic. As a consequence we obtain generators and relations for certain types of fusion products including the $N$--fold fusion product of $V(\lambda)$. This result shows that the fusion product of these types is independent of the chosen parameters, proving a special case of a conjecture by Feigin and Loktev. Moreover, we prove a conjecture by Chari, Fourier and Sagaki on truncated Weyl modules for certain classes of dominant integral weights and show that they are realizable as fusion products. In the last section we consider the case $\mathfrak{g}=\mathfrak{sl}_2$ and compute a PBW type basis for truncated Weyl modules of the associated current algebra.
\end{abstract}
\maketitle \thispagestyle{empty}
%%%%%%%%%%%%%%%%%%%%%%%%%%%%%%%%%%%%%%%%%%%%%%%%%%%%%%%%%%%%%%%%%%%%%%%%%%%%%%%%%%%%%%%%%%%%%%%%%%%%%%%%%%%%%%%%%%%%%%%%%%%%%%%%%%%
%         Introduction
%%%%%%%%%%%%%%%%%%%%%%%%%%%%%%%%%%%%%%%%%%%%%%%%%%%%%%%%%%%%%%%%%%%%%%%%%%%%%%%%%%%%%%%%%%%%%%%%%%%%%%%%%%%%%%%%%%%%%%%%%%%%%%%%%%%
\section{Introduction}

Let $\lie g$ be a finite--dimensional complex simple Lie algebra with highest root $\theta$. The current algebra $\lie g[t]$ associated to $\Lg$ is the algebra of polynomial maps $\mathbb C \rightarrow \lie g$ or equivalently, it is the complex vector space $\Lg\otimes \C[t]$ with Lie bracket the $\C[t]$--bilinear extension of the Lie bracket on $\lie g$. The toroidal current algebra $\lie g[t,u]$ associated to $\Lg$ is the algebra of polynomial maps $\C^2 \rightarrow \lie g$  and can be identified with the complex vector space $\Lg\otimes \C[t,u]$ with similar Lie bracket. The Lie algebra $\Lg[t]$ is graded by the non--negative integers, where the $r$--th graded component is $\Lg\otimes t^r$ and $\lie g[t,u]$ is bi--graded by pairs of non--negative integers, where the $(r,s)$--th graded component is $\Lg\otimes t^ru^s$. We are interested in the category of finite--dimensional graded representations of $\Lg[t]$ and finite--dimensional bi--graded representations of $\lie g[t,u]$ respectively. The former category contains a large number of interesting objects, for example local Weyl modules (see for instance \cite{CIK14,CP01,FoL07,FKKS11}), $\Lg$--stable Demazure modules (see \cite{CSVW14,FoL06,FoL07}) and fusion products.\par The latter class of representations has been introduced by Feigin and Loktev in \cite{FL99}:
given finite--dimensional cyclic $\Lg[t]$--modules $\V_1,\dots,\V_N$ with cyclic vectors $v_1,\dots,v_N$ and a tuple of pairwise distinct complex numbers $\mathbf z=(z_1,\dots,z_N)$ one can define a filtration on the tensor product $\V^{z_1}_1\otimes\cdots \otimes \V^{z_N}_N$ and build the associated graded space with respect to this filtration. This space is called the fusion product and is denoted by $\V^{z_1}_1*\cdots* \V^{z_N}_N$, where $\V^{z}$ is a non--graded $\Lg[t]$--module (see Section~\ref{section3} for more details). The Feigin--Loktev conjecture states that under suitable conditions on $\V_s$ and $v_s$ the fusion product is independent of the chosen fusion parameters $\mathbf z$. This conjecture has been proved for several classes of representations. For example it has been proved in \cite{FoL07,Na11} that the fusion product of local Weyl modules is again a local Weyl module and hence independent of the chosen parameters. Another example are fusion products of Kirillov--Reshetikhin modules (see \cite{Ke11}) and fusion products of $\Lg$--stable Demazure modules (see \cite{CSVW14,FoL07,KV14}).\par
An other interesting class of $\lie g[t]$--modules are those which are obtained as fusion products of finite--dimensional simple $\Lg$--modules, where a $\Lg$--module $\V$ is made into a $\Lg[t]$--module by requiring $\big(\Lg\otimes t\C[t]\big)\V=0$. Hence for any tuple $(\lambda_1,\dots,\lambda_N)$ of dominant integral weights the fusion product $\V^{z_1}(\lambda_1)*\cdots*\V^{z_N}(\lambda_N)$ can be defined and studied.
For these types of representations the Feigin--Loktev conjecture has been proved in the case of $\mathfrak{sl}_2$ and in some special other cases (see for instance \cite{CV13,FF02,FJKLM04,Ra14}). Moreover, in the case of $\mathfrak{sl}_2$ a presentation for the fusion product $\V(k_1)*\cdots*\V(k_N)$ has been established in terms of generators and relations of the enveloping algebra (see \cite{CV13,FF02}) and a PBW type basis has been computed \cite{CV13}. An easy calculation shows that the aforementioned presentation can be greatly simplified if the highest weights are equal. In particular, $\V(k)*\cdots*\V(k)$ is a cyclic $\mathbf{U}(\mathfrak{sl}_2[t])$--module generated by a vector $v$ subject to the same relations as the highest weight vector of the local Weyl module $\W_{\loc}(kN)$ with the only additional relation $\big(\mathfrak{sl}_2\otimes t^N\big) v=0$.\par
This paper is motivated by the idea to generalize the above observation for arbitrary $\lie g$:
is the fusion product $\V^{z_1}(\lambda)*\cdots*\V^{z_N}(\lambda)$ independent of the fusion parameters for arbitrary $\Lg$? Is there a simple presentation considered as a $\mathbf{U}(\Lg[t])$--module? Is the truncated Weyl module $\W(N\lambda,N)$ realizable as a fusion product?
In this paper we give a positive answer to these questions. Our approach is based on the theory of finite--dimensional bi--graded modules for the toroidal current algebra $\Lg[t,u]$. In particular we construct an associated graded version of a $\Lg$--stable Demazure module $\gr_{t^N} \T(\ell,N)$ and a bi--graded version of a fusion product $\overline{\D^{u}(\ell,\ell \lambda)}*\cdots*\overline{\D^{u}(\ell,\ell \lambda)}*\overline{\D^{u}(\ell,\ell \lambda+\lambda^0)}$ such that the zeroth graded space (with respect to the $u$--grading) of the second construction is isomorphic to the fusion product of finite--dimensional simple $\lie g$--modules. Our first result is the following; for the precise definition of the ingredients see Section~\ref{section34}. \textit{We remark that if $\lambda^0\neq 0$, then the Lie algebra $\Lg$ is assumed to be classical or $\tt G_2$, for $\lambda^0=0$ there is no restriction on $\lie g$.}
\begin{theom}Let $\ell\in \mathbb N$, $\lambda\in L^+$ and $\lambda^0\in P^+$ such that $\lambda^0(\theta^{\vee})\leq \ell$. We have an isomorphism of $\mathbf{U}\big(\Lg[t,u]\big)$--modules
$$
\gr_{t^N} \T(\ell,N) \cong \overline{\D^{u}(\ell,\ell \lambda)}*\cdots*\overline{\D^{u}(\ell,\ell \lambda)}*\overline{\D^{u}(\ell,\ell \lambda+\lambda^0)}.
$$
\end{theom}
Our second result gives a connection to truncated Weyl modules, where the first part is a direct consequence of the previous theorem and the second part proves a special case of a conjecture by Chari, Fourier and Sagaki. Again for the precise definition of the ingredients see Section~\ref{section41}. 
\begin{theom}
Let $\ell\in \mathbb N$, $\lambda\in L^+$ and $\lambda^0\in P^+$ such that $\lambda^0(\theta^{\vee})\leq \ell$.
\begin{enumerate}
\item The fusion product $\V(\ell\lambda)^{*(N-1)}*\V(\ell \lambda+\lambda^0)$ is independent of the fusion parameters.
\item If $\lambda^0(\theta^{\vee})\leq 1$ and $|N\lambda+\lambda^0|\geq N$, then 
$$\W(N \lambda+\lambda^0,N)\cong \V(\lambda)^{*(N-1)}*\V(\lambda+\lambda^0).$$
\end{enumerate}
In Section~\ref{section43} and Section~\ref{section44} we deal with the case of $\mathfrak{sl}_2$ and prove that the truncated Weyl module is realizable as a fusion product. Moreover, we compute a PBW type basis for truncated Weyl modules which differs from the basis described in \cite[Section 6]{CV13}. 
\end{theom}
\begin{theom}
 Let $m\in \mathbb Z_{+}$ and write $m=kN+j$,\ $0\leq j<N$.
\begin{enumerate}
\item We have an isomorphism of $\mathbf U\big(\mathfrak{sl}_2\otimes \C[t]/t^N\big)$--modules 
$$\W(m,N)\cong V(k)^{*(N-j)}*V(k+1)^{*j}$$
\item A $\PBW$ type basis of $\W(m,N)$ is given by $$\Big\{(x_{-\alpha}\otimes 1)^{i_0}\cdots (x_{-\alpha}\otimes t^{N-1})^{i_{N-1}}w_{m,N}\mid(i_0,\ldots,i_{N-1})\in S\big(k^{N-j},(k+1)^{j}\big)\Big\}.$$
\end{enumerate}
\end{theom}

Our paper is organized as follows. Section~\ref{section2} establishes the basic notation needed in the rest of the paper. In Section~\ref{section3}, we construct in different ways two bi--graded modules and prove them to be isomorphic. As a consequence we obtain that the fusion product is independent of the choosen parameters.
In Section~\ref{section4}, we give some applications regarding the conjecture on truncated Weyl modules and compute a PBW type basis.
%%%%%%%%%%%%%%%%%%%%%%%%%%%%%%%%%%%%%%%%%%%%%%%%%%%%%%%%%%%%%%%%%%%%%%%%%%%%%%%%%%%%%%%%%%%%%%%%%%%%%%%%%%%%%%%%%%%%%%%%%%%%%%%%%%%
\section{Preliminaries}\label{section2}
\subsection{}\label{section21} Throughout the paper $\mathbb C$ denotes the field of complex numbers, $\mathbb Z$ the ring of integers and $\mathbb Z_+$, $\mathbb N$ the set of non--negative and positive integers respectively. Given any complex Lie algebra $\lie a$ we let $\mathbf U(\lie a)$ be the universal enveloping algebra of $\lie a$. Further, let $\lie a[t]$ be the  Lie algebra of polynomial maps from $\mathbb C$ to $\lie a$ with the obvious pointwise Lie bracket: $$[x\otimes f, y\otimes g]=[x,y]\otimes fg,\ \ x,y\in\lie a,\ \ f,g\in\mathbb C[t].$$ The Lie algebra $\lie a[t]$  and its universal enveloping algebra inherit a grading from the degree grading of $\mathbb C[t]$, thus an element $a_1\otimes t^{r_1}\cdots a_s\otimes t^{r_s}$, $a_j\in\lie a$, $r_j\in\mathbb Z_+$ for $1\le j\le s$ will have grade $r_1+\cdots+r_s$. We shall be interested in $\mathbb Z$--graded modules $V=\oplus_{s\in\mathbb Z}V[s]$ for $\lie a[t]$.

\subsection{}\label{section22} We refer to \cite{K90} for the general theory of affine Lie algebras. 
Let $\lie g$  be a finite--dimensional complex simple Lie algebra and $\widehat{\lie g}$ be the corresponding untwisted affine algebra. We fix $\lie h\subset \widehat{\lie h}$ Cartan subalgebras of $\lie g $ and $\widehat{\lie g}$  respectively and denote by $R$ and $\widehat{R}$ respectively the set of roots of $\lie g$ with respect to $\lie h$, and the set of roots of $\widehat{\lie g}$ with respect to $\widehat{\lie h}$. The corresponding sets of positive and negative roots are denoted as usual by $R^\pm$ and $\widehat R^\pm$ respectively. We fix $\Delta=\{\alpha_1,\dots,\alpha_n\}$ a basis for $R$ such that $\widehat\Delta=\Delta\cup\{\alpha_0\}$ is a basis for $\widehat R$. For $\alpha\in \widehat{R}$, let $\alpha^{\vee}$ be the corresponding coroot.
We fix  $d\in \widehat{\lie h}$ such that $\alpha_0(d)=1$ and $\alpha_i(d)=0$ for $i\neq 0$; $d$ is called the scaling element and it is unique modulo the center of $\widehat{\lie g}$. For $1\le i\leq n$, define $\omega_i\in\lie h^*$ by $\omega_i(\alpha_j^\vee)=\delta_{i,j}$, for $1\leq j\leq n$, where $\delta_{i,j}$ is Kronecker's delta symbol. The element $\omega_i$ is the fundamental weight of $\lie g$ corresponding to $\alpha_i^\vee$. Let $(.,.)$ be the non--degenerate symmetric bilinear form on $\lie h^{*}$ normalized so that the square length of a long root is two. For $\alpha\in R^+$ we set 
$$d_{\alpha}=\frac{2}{(\alpha,\alpha)},\ d_{i}:=d_{\alpha_i}, \mbox{ for $1\leq i \leq n$}.$$
The weight lattice $P$ (resp. $P^+$) is the $\mathbb Z$--span (resp. $\mathbb Z_+$--span) of the fundamental weights. 
The co--weight lattice $L$ is the sublattice of $P$ spanned by the elements $d_i\omega_i$, $1\le i\le n$ and the subset $L^+$ is defined in the obvious way. 
For $\lambda\in P^+$ we define 
$$|\lambda|=\sum^n_{i=1}\lambda(\alpha_i^{\vee})\in \mathbb Z_+.$$

%%%%%%%%%%%%%%%%%%%%%%%%%%%%%%%%%%%%%%%%%%%%%%%%%

\subsection{}\label{section23} 
 Given $\alpha\in \widehat R^+$ let $\widehat{\lie g}_\alpha\subset\widehat{\lie g} $ be the corresponding root space; note that  $\widehat{\lie g}_\alpha \subset \lie g$ if $ \alpha\in R$. For a real root $\alpha$ we denote by $x_{\alpha}$ a generator of $\widehat{\lie g}_\alpha$.
 The element $d$ defines a $\mathbb Z_+$--graded Lie algebra structure on $\Lg[t]$: for $\alpha\in \widehat R$ we say that $\lie g_\alpha$ has grade $k$ if $$[d,x_\alpha]=kx_\alpha$$ 
 or, equivalently, if $\alpha(d)=k$. With respect to this grading, the zero homogeneous component of the current algebra is  $\Lg[t][0]\cong \lie g$ and the subspace spanned by the positive homogeneous components is an ideal denoted by $\Lg[t]_+$. We have a short exact sequence of Lie algebras,
$$0\to \Lg[t]_+\to\Lg[t] \stackrel{\ev_0}{\longrightarrow} \lie g\to 0.$$
Clearly the pull--back of any $\lie g$--module $\V$ by  $\ev_0$ defines the structure of a graded $\Lg[t]$--module on $\V$ and we denote this module by $\ev_0^*\V$.

\subsection{}\label{section24}  For $\lambda\in P^+$, denote  by $\V(\lambda)$ the simple finite--dimensional $\lie g$--module generated by an element $v_\lambda$ with defining relations $$\lie n^+ v_\lambda=0,\ \ \alpha_i^{\vee}v_\lambda=\lambda(\alpha_i^{\vee}) v_\lambda,\ \ (x_{-\alpha_i})^{\lambda(\alpha_i^{\vee})+1} v_\lambda=0,\ \ 1\le i\le n. $$ 
It is well--known that $\V(\lambda)\cong \V(\mu)$ iff $\lambda=\mu$ and that  any finite--dimensional $\lie g$--module is isomorphic to a direct sum of  modules $\V(\lambda)$, $\lambda\in P^+$. If $\V$ is a $\lie h$--semisimple $\lie g$--module (in particular if $\dim \V<\infty$),  we have $$ \V=\bigoplus_{\mu\in\lie h^*}\V_\mu,\ \ \V_\mu=\{v\in \V \mid hv=\mu(h)v,\ \ h\in\lie h\},$$ and we set $\wt \V=\{\mu\in\lie h^*\mid V_\mu\ne 0\}.$ By our previous comments, for any $\lambda\in P^+$ we obtain a graded $\Lg[t]$--module $\ev_0^*\V(\lambda)$.

We also define the local Weyl module $\W_{\loc}(\lambda)$, which is a finite--dimensional $\lie g[t]$--module generated by an element $w_{\lambda}$ with defining reltations 
$$\lie n^+[t] w_\lambda=0,\ \ (\alpha_i^{\vee}\otimes t^s)w_\lambda=\delta_{s,0}\lambda(\alpha_i^{\vee}) w_\lambda,\ \ (x_{-\alpha_i}\otimes 1)^{\lambda(\alpha_i^{\vee})+1} w_\lambda=0,\ \forall s\geq 0, \ 1\le i\le n. $$ 
For more details regarding the theory of local Weyl modules we refer to \cite{CIK14,CL06,CP01,FKKS11,FK11,FoL07,Na11}.

\subsection{}\label{section25} We recall a general construction from \cite{FL99}. Let $\mathbf U(\Lg[t])[k]$ be the homogeneous component of degree $k$ (with respect to the grading induced by $d$) and recall that it is a $\lie g$--module for all $k\in\mathbb{Z}_+$. Suppose now that we are given a  $\Lg[t]$--module $\V$ which is generated by $v$. Define an increasing filtration $0\subset \V^0\subset \V^1\subset\cdots $ of $\lie g$--submodules of $\V$ by $$\V^k=\bigoplus_{s=0}^k \mathbf U(\Lg[t])[s] v.$$ The associated graded vector space $\gr \V$ admits an action of $\Lg[t]$ given by: 
$$
x(v+\V^{k})= xv+ \V^{k+s},\ \ x\in\Lg[t][s],\ \ v\in \V^{k+1}.
$$
Furthermore, $\gr \V$ is a cyclic $\Lg[t]$--module with cyclic generator $\bar v$, the image of $v$ in $\gr \V$. Given a $\Lg[t]$--module $\V$ and $z\in\mathbb \C$, let $\V^z$ be the $\Lg[t]$--module with action 
$$(x\otimes t^r)w=(x\otimes (t+z)^r)w,\ x\in \lie g,\ w\in \V,\ r\in \mathbb Z_+.$$
Starting with finite--dimensional cyclic $\Lg[t]$--modules $\V_1,\dots,\V_N$ with cyclic vectors $v_1,\dots,v_N$ and a tuple of pairwise distinct complex numbers $\mathbf z=(z_1,\dots,z_N)$, the fusion product is defined to be $\V^{z_1}*\cdots* \V^{z_N}:=\gr (\V^{z_1}\otimes\cdots\otimes \V^{z_N})$. It was proved in \cite{FL99} that the tensor product $\V^{z_1}\otimes\cdots\otimes \V^{z_N}$ is cyclic and generated by $v_1\otimes\cdots \otimes v_N$.
Clearly the definiton of the fusion product depends on the parameters $z_s$, $1\le s\le N$. However it is conjectured in \cite{FL99} (and proved in special cases, see \cite{CL06,FF02,FL99,FoL07,
KV14} for example)  that under suitable conditions on $\V_s$ and $v_s$, the fusion product is independent of the choice of the complex numbers. In this paper we cover another class of representations, where the construction of the fusion product is independent of the parameters. To keep the notation as simple as possible we omit almost always the parameters in the notation for the fusion product and write  $\V_1*\cdots*\V_N$ for $\V_1^{z_1}*\cdots*\V_N^{z_N}$.

%%%%%%%%%%%%%%%%%%%%%%%%%%%%%%%%%%%%%%%%%%%%%%%%%%%%%%%%%%%%%%%%%%%%%%%%%%%%%%%%%%%%%%%%%%%%%%%%%%%%%%%%%%%%%%%%%%%%%%%%%%%%%%%%%%%
\section{Filtrations and bi--graded modules}\label{section3}
The aim of this section is to construct two finite--dimensional bi--graded modules in different ways and prove them to be isomorphic. The advantage of this construction is that a comparison of the zeroth graded components leads to a realization of the fusion product associated to rectangular partitions. 
\subsection{}\label{section31}Let us start with our first construction. Let $\lambda\in P^+$ and $\ell\in \mathbb N$. The $\Lg$--stable Demazure module $\D(\ell,\lambda)$ is a finite--dimensional submodule of a level $\ell$ highest weight representation for the affine algebra $\widehat{\lie g}$.
For these representations, generators and relations are known if we consider them as $\mathbf{U}(\Lg[t])$--modules, see \cite{FK11,FoL07,M88} for more details. 
We remark that these relations are greatly simplified for Demazure modules for untwisted affine algebras in \cite{CV13} and for twisted affine algebras in \cite{KV14}. For instance, one can use these simplified relations to see directly that level one Demazure modules are isomorphic to local Weyl modules for simply--laced affine algebras and twisted affine algebras respectively, initially proved in \cite{FoL07} and \cite{CIK14,FK11}  respectively. We recall the simplified presentation of $\Lg$--stable Demazure modules.
Write 
\begin{equation}\label{zer}\lambda(\beta^{\vee})=(p_{\beta}-1)d_{\beta}\ell+m_{\beta},\ 0<m_{\beta}\leq d_{\beta}\ell, \ \mbox{ for $\beta\in R^+$}. \end{equation}

\begin{prop}\label{demrel0}
The Demazure module $\D(\ell,\lambda)$ is isomorphic to the cyclic $\mathbf{U}(\Lg[t])$--module generated by a
vector $v\neq 0$ subject to the following relations:
\begin{align}
\label{e1}&\Ln^{+}[t]v = 0,\quad \big(h\otimes t^{s}\big)v =\delta_{s,0} \lambda(h) v,\ \forall h\in \Lh,\ s\ge 0
&\\&\label{e2}\big(x_{-\beta} \otimes 1\big)^{\lambda(\beta^{\vee})+1}v = 0,\quad \big(x_{-\beta} \otimes t^{p_{\beta}}\big)v = 0,\ \forall \beta\in R^+
&\\&\label{e3}
\big(x_{-\beta} \otimes t^{p_{\beta}-1}\big)^{m_{\beta}+1}v = 0,\ \forall \beta\in R^+ \mbox{ such that $m_{\beta}<d_{\beta}\ell$}.\end{align}
\end{prop}
We can decompose $\D(\ell,\lambda)$ into simple finite--dimensional $\Lg$--modules. We remark that the vector $v$ in Proposition~\ref{demrel0} corresponds to the highest weight vector of $\ev_0^{*}\V(\lambda)$ in the $\Lg$--module decomposition of $\D(\ell,\lambda)$. We call it a highest weight vector of the module.

\subsection{}\label{section32}
In this paper we are concerned with Demazure modules of the form $\D(\ell,\ell N\lambda^1+\lambda^0)$, where $\lambda^1\in L^+$ and $\lambda^0\in P^+$ such that $\lambda^0(\theta^{\vee})\leq \ell$. 
%$\lambda^0(\beta^{\vee})\leq d_{\beta}\ell$ for all $\beta\in R^+$.
\textit{For the rest of this paper we assume that either $\lambda^0=0$ and $\lie g$ is arbitrary or $\lambda^0\neq 0$ and $\lie g$ is of classical type or $\tt G_2$.}
By the results of \cite{CSVW14,FoL07}, the Demazure module $\D(\ell,\ell N\lambda^1+\lambda^0)$ is isomorphic to the fusion product of $(N-1)$--copies of the Demazure module $\D(\ell,\ell\lambda^1)$ with $\D(\ell,\ell\lambda^1+\lambda^0)$:
\begin{equation}\label{fus}
\D(\ell, \ell N\lambda^1+\lambda^0) \cong \D(\ell,\ell\lambda^1) * \cdots * \D(\ell,\ell\lambda^1)* \D(\ell,\ell\lambda^1+\lambda^0).
\end{equation}
This decomposition holds for all fusion parameters $\mathbf{z}=(z_1,\ldots,z_N)$ with $z_i \neq z_j$ for all $i\not=j$. We emphasize that the restriction on $\lie g$ is only made because \eqref{fus} is not proved for the remaining exceptional Lie algebras if $\lambda^0$ is non--zero. In other words, our results are applicable whenever we have such a fusion product decomposition. We will need the following lemma.

\begin{lem}\label{lem2}
Let $\beta$ be a positive root and $\lambda\in P^+$. We write $(\theta-\beta)=\sum_j \gamma_j$ as a sum of positive roots. Then we have
$$
\lambda(\beta^{\vee})(\beta,\beta)\leq \lambda(\theta^{\vee})(\theta,\theta) \text{ and equality holds iff}\ \lambda\big(\gamma_j^{\vee}\big)=0 \ \forall j.
$$
\proof
Since $\lambda$ is a dominant integral weight we have
$\lambda(\beta^{\vee})\ge 0$ for a positive root $\beta$. We obtain
$$\theta^{\vee}=(\beta+\sum_j\gamma_j)^{\vee}=\frac{(\beta,\beta)}{2}\beta^{\vee}+\sum_j\frac{(\gamma_j,\gamma_j)}{2}\gamma_j^{\vee},$$
which gives $$\lambda(\theta^{\vee})(\theta,\theta) =(\beta,\beta)\lambda(\beta^{\vee})+\sum_j(\gamma_j,\gamma_j) \lambda(\gamma_j^{\vee})\geq (\beta,\beta)\lambda(\beta^{\vee}).$$
Note that equality is only possible if $\lambda(\gamma_j^{\vee})=0$ for all $j$, since $(\gamma_j,\gamma_j)>0$.
\endproof
\end{lem}
By Lemma~\ref{lem2} and Proposition~\ref{demrel0}, equation~\eqref{e2} we get:
\begin{cor}
$\big(x_{-\beta} \otimes t^{(\lambda^1(\theta^{\vee})+1)N}\big)v=0$ for all roots $\beta\in R^+$.
\proof
Write $(\ell N\lambda^1+\lambda^0)(\beta^{\vee})$ as in \eqref{zer}. Since $\lambda^1\in L^+$ we have $m_{\beta}=d_{\beta}\ell$ if $\lambda^0(\beta^{\vee})=0$ and $m_{\beta}=\lambda^0(\beta^{\vee})$ else. Then $\big( x_{-\beta} \otimes t^{p_{\beta}}\big)v=0$ and 
$$p_{\beta}=N\frac{\lambda^1(\beta^{\vee})}{d_{\beta}}+\frac{\lambda^0(\beta^{\vee})-m_{\beta}}{d_{\beta}\ell}+1\leq N(\lambda^1(\theta^{\vee})+1).$$
\endproof
\end{cor} 
Hence $\D(\ell,\ell N\lambda^1+\lambda^0)$ is a $\mathbf{U}\big(\Lg\otimes \C[t]/t^{(\lambda^1(\theta^{\vee})+1)N}\big)$--module.
\subsection{}\label{section33}
We define a decreasing filtration on $\mathbf{U}\big(\Lg\otimes \C[t]/t^{(\lambda^1(\theta^{\vee})+1)N}\big)$
$$ \T_0(N) \supseteq \T_1(N) \supseteq \T_2(N) \supseteq \cdots$$
with 
$$
\T_0(N)=\mathbf{U}\big(\Lg\otimes \C[t]/t^{(\lambda^1(\theta^{\vee})+1)N}\big),\ 
\text{and}\ \T_j(N)=\big(\Lg\otimes t^N\C[t]\big)\T_{j-1}(N) \ \text{for\ }j\ge 1
$$
and study the induced decreasing filtration on our Demazure module $\D(\ell,\ell N\lambda^1+\lambda^0)$ given by
$$
\D(\ell,\ell N\lambda^1+\lambda^0) = \T_0(N)v=:\T_0(\ell,N) \supseteq \T_1(N)v:=\T_1(\ell,N) \supseteq \T_2(\ell,N) \supseteq \cdots
$$
To be consistent with the notation in \cite{F2008}, we refer to it as the $t^N$--filtration.
Let $\gr_{t^N}\T(N)$ respectively $\gr_{t^N} \T(\ell,N)$ be the associated graded space:
\begin{align*}
\gr_{t^N} T(N)\hspace{0,65cm} &= \T_0(N)/\T_1(N) \oplus \T_1(N)/\T_2(N) \oplus \ldots \\
\text{respectively} \hspace{0,2cm}\gr_{t^N} \T(\ell,N) &= \T_0 (\ell,N)/\T_1(\ell,N) \oplus \T_1(\ell,N)/\T_2(\ell,N) \oplus \ldots 
\end{align*}
Since $\D(\ell,\ell N\lambda^1+\lambda^0)$ is a module for $\mathbf{U}\big(\Lg\otimes \C[t]/t^{(\lambda^1(\theta^{\vee})+1)N}\big)$ we obtain that $\gr_{t^N} \T(\ell,N)$ is a module for 
$\gr_{t^N} \T(N)$. By the following lemma $\gr_{t^N} \T(\ell,N)$ is also a module for the toroidal current algebra $\mathbf{U}\big(\Lg\otimes \C[t,u]/\langle t^N,u^{\lambda^1(\theta^{\vee})+1}\rangle\big)$. 
\begin{lem}\label{isoal}
We have an isomorphism of algebras
$$\Psi:\gr_{t^N} \T(N)\stackrel{\sim}{\longrightarrow} \mathbf{U}\big(\Lg\otimes \C[t,u]/\langle t^N,u^{\lambda^1(\theta^{\vee})+1}\rangle\big),$$
where $\Psi\big(x\otimes t^{jN+s}\big)=x\otimes u^jt^s$ for $x\in\Lg$ and  $0\leq s<N$.
\proof
The map $\Psi$ is clearly an isomorphism of vector spaces. In order to show that this map is an algebra homomorphism, we have to check that the naive way of defining $\Psi$ on a product of elements is well--defined. Hence we will verify that 
$$
\big(x\otimes u^jt^s\big)\big(y\otimes u^it^q\big)-\big(y\otimes u^it^q\big)\big(x\otimes u^jt^s\big)=\big[x,y\big]\otimes u^{i+j}t^{s+q}
$$
holds on the right hand side, whenever we have 
$$
\big(x\otimes t^{jN+s}\big)\big(y\otimes t^{iN+q}\big)-\big(y\otimes t^{iN+q}\big)\big(x\otimes t^{jN+s}\big)=\big[x,y\big]\otimes t^{(i+j)N+(s+q)}
$$
on the left hand side. This is obvious for $(s+q)<N$. Otherwise the variables $x\otimes u^jt^s$ and $y\otimes u^it^q$ commute in $\mathbf{U}\big(\Lg\otimes \C[t,u]/\langle t^N,u^{\lambda^1(\theta^{\vee})+1}\rangle\big)$. By the definition of the associated graded space we also obtain that the variables $x\otimes t^{jN+s}$ and $y\otimes t^{iN+q}$ commute in $\gr_{t^N} \T(N)$ since on the one hand 
$$\big(x\otimes t^{jN+s}\big)\big(y\otimes t^{iN+q}\big)-\big(y\otimes t^{iN+q}\big)\big(x\otimes t^{jN+s}\big)\in \T_{i+j}(N)$$ and on the other hand $$\big[x,y\big]\otimes t^{(i+j)N+(s+q)}\in \T_{i+j+1}(N).$$
\endproof
\end{lem}
\subsection{}\label{section34}
Now we present a quite different construction of the module $\gr_{t^N} \T(\ell,N)$. In fact, it is one of the main results of this paper to prove that the two constructions give isomorphic modules. We start with the $(N-1)$--fold tensor product of Demazure modules $\D(\ell,\ell\lambda^1)$ with $\D(\ell,\ell\lambda^1+\lambda^0)$. The gambit: we switch the variables and consider now the current algebra $\lie g[u]$ which operates on the Demazure modules $\D^u(\ell,\ell\lambda^1)$. Here we add the index $u$ to emphazise that here the algebra $\lie g[u]$ is acting. We extend the action trivially to $\lie g[t,u]$ and denote the corresponding module by $\overline{\D^u(\ell,\ell\lambda^1)}$, i.e. $\lie g\otimes t\C[t,u]$ acts trivially. Recall that we get a highly non-trivial action of $\lie g[t,u]$ when we consider fusion products of these modules with respect to the variable $t$. The bi-graded fusion product 
\begin{equation}\label{bigrfus}\overline{\D^{u}(\ell,\ell \lambda^1)}*\cdots*\overline{\D^{u}(\ell,\ell \lambda^1)}*\overline{\D^{u}(\ell,\ell \lambda^1+\lambda^0)}\end{equation} is a cyclic module for the Lie algebra $\mathbf U(\Lg \otimes \C[t,u]/\langle t^N, u^{\lambda^1(\theta^{\vee})+1}\rangle)$. 
Note the similarity but also the difference between \eqref{fus} and \eqref{bigrfus}. In \eqref{fus} we consider the fusion product (with respect to the variable $t$) of $\lie g[t]$--modules. The $\lie g[t,u]$--module structure comes only into the picture by the filtration defined in Section~\ref{section33}. We would like to remind the reader that if $\lambda^0\neq 0$, then $\Lg$ is of classical type or $\tt G_2$.
\begin{thm}\label{1}
Let $\lambda^1\in L^+$ and $\lambda^0\in P^+$ such that $\lambda^0(\theta^{\vee})\leq \ell$. We have an isomorphism of $\mathbf{U}\big(\Lg \otimes \C[t,u]/\langle t^N, u^{\lambda^1(\theta^{\vee})+1}\rangle\big)$--modules
$$
\gr_{t^N} \T(\ell,N) \cong \overline{\D^{u}(\ell,\ell \lambda^1)}*\cdots*\overline{\D^{u}(\ell,\ell \lambda^1)}*\overline{\D^{u}(\ell,\ell \lambda^1+\lambda^0)}.
$$
\proof
Let $v_{\ell}^{* (N-1)}*v_0$ be the highest weight vector of the right hand side. The isomorphism between $\gr_{t^N} \T(N)$ and $\mathbf{U}\big(\Lg \otimes \C[t,u]/\langle t^N, u^{\lambda^1(\theta^{\vee})+1}\rangle\big)$ (see Lemma~\ref{isoal}) induces a natural surjective map 
$$
\gr_{t^N} \T(N)\twoheadrightarrow \mathbf{U}\big(\Lg \otimes \C[t,u]/\langle t^N, u^{\lambda^1(\theta^{\vee})+1}\rangle\big)\circ (v_{\ell}^{*(N-1)}*v_0).
%=\overline{\D^{u}(\ell,\ell \lambda^1)}*\cdots*\overline{\D^{u}(\ell,\ell \lambda^1)}*\overline{\D^{u}(\ell,\ell \lambda^1+\lambda^0)}.
$$
It remains to prove that this map induces an isomorphism between the cyclic module
$\gr_{t^N} \T(\ell,N)$ and the fusion product.
Since the dimensions of the modules coincide it is enough to show that all relations which hold in 
$\gr_{t^N} \T(\ell,N)$ also hold on the right hand side.\par
Recall that a presentation of $\gr_{t^N} \T(\ell,N)$ is given by two types of relations, the ones coming from the presentation of the Demazure module and the ones coming from going to the associated graded space with respect to the $t^N$--filtration.
We start by proving that the defining relations of $\D(\ell,\ell N\lambda^1+\lambda^0)$ given for $v$ in Theorem~\ref{1} are satisfied by $v_{\ell}*\cdots*v_{\ell}*v_0$.
Since the relations \eqref{e1} and the first part of \eqref{e2} are obviously satisfied it remains to verify the second part of \eqref{e2} and \eqref{e3}. Write $(\ell N\lambda^1+\lambda^0)(\beta^{\vee})$ as in \eqref{zer}. We start to prove
\begin{equation}\label{demrel}\Big(x_{-\beta}\otimes u^{j_{\beta}}t^{r_{\beta}}\Big)(v_{\ell}^{* (N-1)}*v_0)=0,\ \mbox{where $p_{\beta}=j_{\beta}N+r_{\beta}$, $0\leq r_{\beta} <N$}.\end{equation}
Since $\lambda^0(\beta^{\vee})\leq d_{\beta}\ell$, we have $$p_{\beta}=\begin{cases}N\lambda^1(\beta^{\vee})d^{-1}_{\beta},& \text{if $\lambda^0(\beta^{\vee})=0$}\\
N\lambda^1(\beta^{\vee})d^{-1}_{\beta}+1,& \text{else}.\end{cases}$$
In either case $j_{\beta}\geq \lambda^1(\beta^{\vee})d^{-1}_{\beta}$ and thus $\big(x_{-\beta}\otimes u^{j_{\beta}}t^{r_{\beta}}\big)v_{\ell}=0$ follows immediately from the defining relations of $\D(\ell,\ell \lambda^1)$. If $r_{\beta}\neq 0$ we can replace $t^{r_{\beta}}$ by $(t-z_N)^{r_{\beta}}$ in the associated graded space and obtain that the element in \eqref{demrel} acts trivially on $v_0$. If $r_{\beta}=0$, then $p_{\beta}$ is divisible by $N$ which forces $\lambda^0(\beta^{\vee})=0$. Therefore, in this case we obtain $j_{\beta}= \lambda^1(\beta^{\vee})d^{-1}_{\beta}$ and $(x_{-\beta}\otimes u^{j_{\beta}})v_0=0$ follows immediately from the defining relations of $\D(\ell,\ell \lambda^1+\lambda^0)$. It remains to consider the relations \eqref{e3}. So suppose we have 
$$p_{\beta}-1=N\frac{\lambda^1(\beta^{\vee})}{d_{\beta}}+\frac{\lambda^0(\beta^{\vee})-m_{\beta}}{d_{\beta}\ell}=j_{\beta}N+r_{\beta},\ 0\leq r_{\beta} <N.$$
Since $m_{\beta}< d_{\beta}\ell$, we must have $m_{\beta}=\lambda^0(\beta^{\vee})\neq 0$ and hence $p_{\beta}-1=N\lambda^1(\beta^{\vee})d^{-1}_{\beta}$.
It follows $j_{\beta}=\lambda^1(\beta^{\vee})d^{-1}_{\beta}$ and therefore $(x_{-\beta}\otimes u^{j_{\beta}})v_{\ell}=0$. It means that $$\Big(x_{-\beta}\otimes u^{j_{\beta}}\Big)^{m_{\beta}+1}(v_{\ell}^{* (N-1)}*v_0)=v_{\ell}^{* (N-1)}*\Big(x_{-\beta}\otimes u^{j_{\beta}}\Big)^{m_{\beta}+1}v_0=0.$$
We consider now the relations coming from the $t^N$--filtration. Suppose 
$$
M=\sum_m\sum_{\substack{i_1,\dots,i_m\\j_1,\dots,j_m}} k(m)_{\substack{i_1,\dots,i_m\\j_1\dots,j_m}} \big(x_{i_1}\otimes t^{i_1N+j_1}\big)
\cdots \big(x_{i_m}\otimes t^{i_mN+j_m}\big)\in \mathbf{U}\big(\Lg\otimes \C[t]/t^{(\lambda(\theta^{\vee})+1)N}\big)
$$
is a linear combination of monomials with fixed $t^N$--degree such that $w=Mv\not=0$ in $\D(\ell,\ell N\lambda^1+\lambda^0)$ but
the image $\overline{w}=0$ in $\gr_{t^N} \T(\ell,N)$. This is only possible if there exists a linear combination of monomials of greater $t^N$--degree 
$$
M'=\sum_{m'}\sum_{\substack{p_1,\dots,p_{m'}\\q_1\dots,q_{m'}}}
k(m^{'})_{\substack{p_1,\dots,p_{m'}\\q_1,\dots,q_{m'}}} \big(x_{p_1}\otimes t^{p_1N+q_1}\big)
\cdots \big(x_{p_{m'}}\otimes t^{p_{m'}N+q_{m'}}\big)\in \mathbf{U}\big(\Lg\otimes \C[t]/t^{(\lambda(\theta^{\vee})+1)N}\big)
$$
such that $w=Mv=M'v$ in $\D(\ell,\ell N\lambda^1+\lambda^0)$. We assume in the following that $M'$ is of maximal $t^N$--degree.
We have $(M-M')v=0$, so the difference $(M-M')$ is an element in the left ideal generated by the elements in \eqref{e1}--\eqref{e3}.
Since $M'$ is of higher $t^N$-degree we get in $\gr_{t^N} \T(N)$:
$$
\overline{(M-M')}=\overline{M},
$$
and since all defining relations of $\D(\ell,\ell N\lambda^1+\lambda^0)$ are satisfied by $(v_{\ell}^{* (N-1)}*v_0)$ we get: $\Psi(\overline{M} )\circ (v_{\ell}^{* (N-1)}*v_0)=0,$
which shows that 

the natural surjective map 
$$
\gr_{t^N} \T(N)\twoheadrightarrow \mathbf{U}\big(\Lg \otimes \C[t,u]/\langle t^N, u^{\lambda^1(\theta^{\vee})+1}\rangle\big)\circ (v_{\ell}^{* (N-1)}*v_0)
$$
induces an isomorphism of cyclic modules $\gr_{t^N} \T(\ell,N) \cong \overline{\D^{u}(\ell,\ell \lambda^1)}*\cdots*\overline{\D^{u}(\ell,\ell \lambda^1)}*\overline{\D^{u}(\ell,\ell \lambda^1+\lambda^0)}$.
\endproof
\end{thm}
For the rest of this section we discuss a crucial consequence of our result.
\begin{cor}\label{inde}
Let $\ell\in \mathbb N$, $\lambda^1\in L^+$ and $\lambda^0\in P^+$ such that $\lambda^0(\theta^{\vee})\leq \ell$. 
\begin{enumerate}
\item The fusion product $\V(\ell\lambda^1)^{*(N-1)}*\V(\ell \lambda^1+\lambda^0)$ is independent of the fusion parameters $\mathbf{z}$.
\item We have an isomorphism of $\mathbf U\big(\Lg\otimes \C[t]/t^N\big)$--modules
$$\V(\ell\lambda^1)^{*(N-1)}*\V(\ell \lambda^1+\lambda^0)\cong \D(\ell,\ell N\lambda^1+\lambda^0)/\big(\Lg\otimes t^N\C[t]\big)\D(\ell,\ell N\lambda^1+\lambda^0).$$
\item If $\lambda^0(\theta^{\vee})\leq 1$, the truncated level one Demazure module is isomorphic to the truncated level $\ell$ Demazure module
$$\D(1,\ell N\lambda^1+\lambda^0)/\big(\Lg\otimes t^N\C[t]\big)\D(1,\ell N\lambda^1+\lambda^0)\cong\D(\ell,\ell N\lambda^1+\lambda^0)/\big(\Lg\otimes t^N\C[t]\big)\D(\ell,\ell N\lambda^1+\lambda^0).$$
\end{enumerate}
\proof
Since the fusion product $\V(\ell\lambda^1)^{*(N-1)}*V(\ell \lambda^1+\lambda^0)$ is isomorphic to the zeroth graded component of $\overline{\D^{u}(\ell,\ell \lambda^1)}*\cdots*\overline{\D^{u}(\ell,\ell \lambda^1)}*\overline{\D^{u}(\ell,\ell \lambda^1+\lambda^0)}$ (with respect to the $u$--grading) the statement follows from Theorem~\ref{1}.
\endproof
\end{cor}
\begin{rem} 
Theorem~\ref{1} generalizes a result of \cite{F2008}, where the statement of the theorem is proved for $\ell=1$, $\lambda^0=0$ and $\lambda^1=\theta$. Unfortunately, the proof in \cite{F2008} has a gap (personal communication by the author), which is now fixed by the proof above. The author of \cite{Ra14} uses the result of \cite{F2008} to prove a presentation for the fusion product $\V(\theta)^{*N}*\D(1,\theta)^{*M}$. 
\end{rem}
%%%%%%%%%%%%%%%%%%%%%%%%%%%%%%%%%%%%%%%%%%%%%%%%%%%%%%%%%%%%%%%%%%%%%%%%%%%%%%%%%%%%%%%%%%%%%%%%%%%%%%%%%%%%%%%%%%%%%%%%%%%%%%%%%%%
%         
%%%%%%%%%%%%%%%%%%%%%%%%%%%%%%%%%%%%%%%%%%%%%%%%%%%%%%%%%%%%%%%%%%%%%%%%%%%%%%%%%%%%%%%%%%%%%%%%%%%%%%%%%%%%%%%%%%%%%%%%%%%%%%%%%%%
\section{Truncated Weyl modules and PBW type basis}\label{section4}
In this section we shall give some evidence for the conjecture made by Chari, Fourier and Sagaki on truncated Weyl modules. This conjecture is nowhere written in the literature, so we decided so state the conjecture in this paper. Finally, we consider the case $\lie g=\mathfrak{sl}_2$ and compute a PBW type basis.
\subsection{}\label{section41}
Let $P^+(\lambda,N)$ be the set of $N$--tuples of dominant integral weights $\blambda=(\lambda_1,\dots,\lambda_N)$, such that $\sum_i \lambda_i=\lambda$. Let $\blambda=(\lambda_1,\dots,\lambda_N),\bmu=(\mu_1,\dots,\mu_N)\in P^+(\lambda,N)$. For a positive root $\beta$ define
$$r_{\beta,k}(\blambda)=\min\big\{(\lambda_{i_1}+\cdots+\lambda_{i_k})(\beta^{\vee})\mid 1\leq i_1<\cdots<i_k\leq N\big\}.$$
We say $\blambda \preceq \bmu$ if
$$r_{\beta,k}(\blambda)\leq r_{\beta,k}(\bmu) \mbox{ for all $\beta\in R^{+}$ and } 1\leq k\leq N.$$ 
The above partial order was considered in \cite{CFS12}, where the authors made the following observation: for a tuple $\blambda$ the dimension of the tensor product of the corresponding finite--dimensional simple $\Lg$--modules increases along $\preceq$.
Moreover, they prove in certain cases (for instance when $\lambda$ is a multiple of a fundamental minuscule weight) that there exists an inclusion of tensor products along with the partial order and conjecture that this remains true for $N=2$ and arbitrary $\lambda$ (see \cite[Conjecture 2.3]{CFS12}). Using the unique maximal element in the partial ordered set $P^+(\lambda,N)$ one can formulate a conjecture on truncated Weyl modules, which we will explain now. 
\begin{defn}
Let $\lambda\in P^+$. The truncated Weyl module $\W(\lambda,N)$ is a cyclic module for $\mathbf U(\lie g \otimes \C[t]/t^N)$ generated by $w_{\lambda,N}$ with relations:
\begin{align}
\label{eq1}&\big(\Ln^{+} \otimes \C[t]/t^N\big)w_{\lambda,N} = 0,\quad \big(h\otimes t^{s}\big)w_{\lambda,N} =\delta_{s,0} \lambda(h) w_{\lambda,N},\ \forall h\in \Lh,\ s\ge 0
&\\&\label{eq2}\big(x_{-\beta} \otimes 1\big)^{\lambda(\beta^{\vee})+1}w_{\lambda,N} = 0,\ \forall \beta\in R^+.
\end{align}
\end{defn}
The following conjecture gives a connection between truncated Weyl modules and fusion products of irreducible finite--dimensional $\lie g$--modules.
\begin{conj}
Let $\lambda\in P^+$, such that $|\lambda|\geq N$ and $\blambda=(\lambda_1,\dots,\lambda_N)$ be the unique maximal element in $P^+(\lambda,N)$. Then we have an isomorphism of $\mathbf U\big(\lie g \otimes \C[t]/t^N\big)$--modules
$$\W(\lambda,N)\cong V(\lambda_1)*\cdots*V(\lambda_N).$$
\end{conj}
The following result proves the above conjecture for certain classes of dominant integral weights. 
\begin{thm}\label{2}
Let $\lambda\in L^+$ and $\lambda^0\in P^+$ such that $\lambda^0(\theta^{\vee})\leq 1$ and $|N\lambda+\lambda^0|\geq N$. Then we have an ismorphism of $\mathbf U\big(\lie g \otimes \C[t]/t^N\big)$--modules
$$\W(N\lambda+\lambda^0,N)\cong V(\lambda)*\cdots*V( \lambda)*V( \lambda+\lambda^0).$$
\proof
If $\lambda=0$, there is nothing to prove. By Corollary~\ref{inde} we obtain that $$V( \lambda)*\cdots*V(\lambda)*V(\lambda+\lambda^0)\cong \D(1, N\lambda+\lambda^0)/\big(\Lg\otimes t^N\C[t]\big)\D(1, N\lambda+\lambda^0).$$ We have to show that the defining relations of $\D(1,N\lambda+\lambda^0)/\big(\Lg\otimes t^N\C[t]\big)\D(1, N\lambda+\lambda^0)$ hold in the truncated Weyl module. We shall prove only the non obvious relations. Let $\beta\in R^+$ and write $( N\lambda+\lambda^0)(\beta^{\vee})$ as in \eqref{zer}. Then, as before
$$p_{\beta}-1=\begin{cases}
N\lambda(\beta^{\vee})d_{\beta}^{-1},& \text{if $\lambda^0(\beta^{\vee})\neq 0$}\\
N\lambda(\beta^{\vee})d_{\beta}^{-1}-1,& \text{else.}\end{cases}$$
We consider four cases. If $\lambda(\beta^{\vee})\neq 0$ and $\lambda^0(\beta^{\vee})\neq 0$, then $p_{\beta}\geq p_{\beta}-1\geq N$ and hence $$\big(x_{-\beta}\otimes t^{p_{\beta}}\big)w_{ N\lambda+\lambda^0,N}=\big(x_{-\beta}\otimes t^{p_{\beta}-1}\big)w_{ N\lambda+\lambda^0,N}=0.$$ If $\lambda(\beta^{\vee})\neq 0$ and $\lambda^0(\beta^{\vee})=0$, then $p_{\beta}\geq N$ and $m_{\beta}=d_{\beta}$ (recall that \eqref{e3} was only considered when $m_{\beta}< d_{\beta}$). If $\lambda(\beta^{\vee})=0$ and $\lambda^0(\beta^{\vee})=0$, there is nothing to show; so consider the last case $\lambda(\beta^{\vee})=0$ and $\lambda^0(\beta^{\vee})\neq 0$. In this case $p_{\beta}=1$ and $m_{\beta}=\lambda^0(\beta^{\vee})$. Thus we have to prove
$$\big(x_{-\beta}\otimes t\big)w_{ N\lambda+\lambda^0,N}=\big(x_{-\beta}\otimes 1\big)^{m_{\beta}+1}w_{ N\lambda+\lambda^0,N}=0,$$
where the last equation is clear. Note that it is enough to prove that $\big(x_{-\beta}\otimes t\big)$ acts by zero on the highest weight vector of the local Weyl module $\W_{\loc}(N \lambda+\lambda^0)$. Since $\W_{\loc}(N \lambda+\lambda^0)\cong \W^{z_1}_{\loc}(N\lambda)*\W^{z_2}_{\loc}(\lambda^0)$ we get 
$$\big(x_{-\beta}\otimes t\big)(w_{N\lambda}*w_{\lambda^0})= \big(x_{-\beta}\otimes (t-z_2)\big)(w_{ N\lambda}*w_{\lambda^0})=w_{ N\lambda}*\big(x_{-\beta}\otimes t\big)w_{\lambda^0}.$$ If $\Lg$ is not of type $\tt G_2$, then $\W_{\loc}(\lambda^0)$ is irreducible and the statement follows. If $\Lg$ is $\tt G_2$ it is easy to see that the only positive root $\beta$ with $\big(x_{-\beta}\otimes t\big)w_{\lambda^0}\neq 0$ is the longest short root $\beta=\alpha_1+2\alpha_2$. But then $\lambda(\beta^{\vee})\neq 0$.
\endproof
\end{thm}
We shall show that $(\lambda,\dots,\lambda, \lambda+\lambda^0)$ is in fact the unique maximal element in $P^+(N\lambda+\lambda^0,N)$. Since $\lambda^{0}(\theta^{\vee})\leq 1$, there exists at most one simple root $\alpha$ such that $\lambda^0(\alpha^{\vee})>0$. Without loss of generality we suppose $\lambda^0(\alpha_j^{\vee})=0$ for all $j>1$. Assume that $(\mu_1,\dots,\mu_N)\in P^+(N\lambda+\lambda^0,N)$ such that $(\lambda,\dots,\lambda, \lambda+\lambda^0)\preceq (\mu_1,\dots,\mu_N)$. We fix a simple root $\alpha_j$ and a permutation $\sigma_{j}$ such that $$\mu_{\sigma_{j}(1)}(\alpha_j^{\vee})\leq \cdots \leq \mu_{\sigma_{j}(N)}(\alpha_j^{\vee}).$$ We write $\mu_{\sigma_{j}(i)}(\alpha_j^{\vee})=\lambda(\alpha_j^{\vee})+\epsilon_i(j)+\delta_{i,N}\lambda^0(\alpha_j^{\vee})$ for integers $\epsilon_i(j)$. By our assumptions we obtain
$$0\leq \epsilon_1(j)\leq \cdots \leq \epsilon_{N-1}(j) \leq \epsilon_N(j)+\lambda^0(\alpha_j^{\vee})\ \mbox{ and } \ \sum_{p=1}^N\epsilon_p(j)=0.$$
Hence, up to a permutation we have $\mu_i=\lambda$ for $1\leq i\leq N-1$ and $\mu_N=\lambda+\lambda^0$.

\subsection{}\label{section42}
For the rest of this section we prove the conjecture for $\mathfrak{sl}_2$ and compute a PBW type basis. For $0\leq j < N$, let $S\big(k^{N-j},(k+1)^{j}\big)$ be the set of tuples $(i_0,\ldots,i_{N-1})$ satisfying \begin{equation}\label{defgl1}\sum^{N-1}_{p=0}\frac{N!}{N-p}i_p\leq N!k-\sum^{N-4}_{\ell=0}\frac{N!}{(N-\ell)!}(N-\ell-2)!b_\ell+j(N-1)!\end{equation} for integers $b_\ell$ defined as follows:
$0\le b_\ell< N-\ell$ and
$$i_0-j=b_0\hspace{-0,3cm}\mod N,\ i_\ell+(b_{\ell-1}\hspace{-0,3cm}\mod N-\ell)=b_\ell\hspace{-0,3cm} \mod N-\ell,\ \mbox{ for } \ell=1,\dots,N-4.$$ 
The theorem we shall prove is the follwing:
\begin{thm}\label{3}
Let $m\in \mathbb Z_{+}$ and write $m=kN+j$ for $0\leq j <N$. 
\begin{enumerate}
\item We have an isomorphism of $\mathbf U\big(\mathfrak{sl}_2\otimes \C[t]/t^N\big)$--modules
$$\W(m,N)\cong V(k)^{*(N-j)}*V(k+1)^{*j}.$$ 
\item A $\PBW$ type basis of $\W(m,N)$ is given by $$\Big\{(x_{-\alpha}\otimes 1)^{i_0}\cdots (x_{-\alpha}\otimes t^{N-1})^{i_{N-1}}w_{m,N}\mid(i_0,\ldots,i_{N-1})\in S\big(k^{N-j},(k+1)^{j}\big)\Big\}.$$
\end{enumerate}
Similar as above, a simple calculation shows that $(k,\dots,k,k+1\dots,k+1)\in P^+(m,N)$ is in fact the unique maximal element.
\end{thm}
The rest of this section is dedicated to the proof of Theorem~\ref{3}.
\subsection{}\label{section43}
We start to prove the first part of the theorem. A presentation of the fusion product as a $\mathbf{U}(\mathfrak{sl}_2\otimes \C[t])$ has been presented in \cite{CV13}. So by their results it is enough to show that the highest weight vector of $\W(m,N)$ satisfies the defining relations of $V(k)^{*(N-j)}*V(k+1)^{*j}$ given in \cite[Proposition 2.7]{CV13} which are
$$\mathbf{x}_{-\alpha}(r,s)=\sum_{(b_p)_{p\geq 0}\in \mathbf S(r,s)} (x_{-\alpha}\otimes 1)^{b_{0}}(x_{-\alpha}\otimes t)^{b_{1}}\cdots (x_{-\alpha}\otimes t^{s})^{b_{s}},\ s,r,\ell\in \mathbb N,\ r+s\geq 1+r\ell+q+p,$$
where $q=\max\{0,(N-\ell)k\}$, $p=\max\{0,j-\ell\}$ and $\mathbf S(r,s)$ is the set of tuples $(b_p)_{p\geq 0}$ satisfying $b_0+\cdots+b_s=r$ and $b_1+2b_2+\cdots+sb_{s}=s$.
We assume that $r+s\leq m$, because otherwise the claim follows from the following result of Garland \cite{G78}:
$$(x_{\alpha}\otimes t)^{(s)}(x_{-\alpha}\otimes 1)^{(s+r)}-(-1)^{s} \mathbf{x}_{-\alpha}(r,s)\in\mathbf U(\Lg[t])\Ln^+[t].$$  
Our aim is to prove that for any tuple $(b_p)_{p\geq 0}\in \mathbf S(r,s)$ there exists $p\geq N$ such that $b_p\neq 0$. Assume this is not the case. If $\ell\geq N$ we obtain
$$rN\geq r+s \geq 1+r\ell\geq 1+rN,$$ which is obviously a contradiction. So assume $l\leq N-1$. It follows
$$m\geq r+s \geq 1+r\ell+(N-\ell)k+p=1+\ell(r-k)+m-j+p$$
and thus $r\leq k$. Therefore we obtain the following contradiction
$$1+\ell(r-k)+m-j+p\leq r+s\leq rN\Rightarrow 1\leq (N-\ell)(r-k)-p.$$
Hence $$\W(m,N)\cong V(k)^{*(N-j)}*V(k+1)^{*j}.$$

\subsection{}\label{section44}
Now we will prove the second part of the theorem. For simplicity we write $f_i$ for $x_{-\alpha}\otimes t^i,\ 1\leq i\leq N-1$ and consider the map $\sh:\mathbf{U}(\Ln^{-}[t])\rightarrow \mathbf{U}(\Ln^{-}[t])$ given by $\sh(f_i)=f_{i+1}$. We will need the following result from \cite{FF02}.
\begin{prop}\label{fp}
Let $k_1\leq k_2\leq \cdots \leq k_N$. We have a short exact sequence of $\mathbf{U}(\Ln^{-}[t])$--modules 
$$0\rightarrow V(k_1)*\cdots *V(k_{N-1})\stackrel{\sh}{\longrightarrow} V(k_1)*\cdots *V(k_N)\stackrel{f_0^{-1}}{\longrightarrow} V(k_1)*\cdots *V(k_N-1)\rightarrow 0.$$
\end{prop}
Due to Proposition~\ref{fp} one can construct inductively a $\PBW$ type basis of the fusion product. To be more precise, we have
\begin{equation}\label{un}B(k_1,\ldots,k_N)=B(k_1,\ldots,k_{N-1})_{\sh}\cup f_0 B(k_1,\ldots,k_N-1),\end{equation} where $B(-)$ denotes a basis of the appropriate fusion product.
\begin{example}We have
$$B(1,2)=B(1)_{\sh}\cup f_0B(1,1)$$ 
and hence
$$B(1,2)=\{1,f_0\}_{\sh}\cup f_0\{1,f_0,f_0^2,f_1\}=\{1,f_1,f_0f_1, f_0,f_0^2,f_0^3\}.$$
\end{example}
We have the following recursion formula:
\begin{lem}\label{mainlem}
We obtain
\begin{align*}
B(k^N)=\bigcup^{k}_{r=0}f^{Nr}_0 B\big((k-r)^{N-1}\big)_{\sh}\cup \bigcup^{N-1}_{j=1}\bigcup^{k}_{r=1} f^{Nr-j}_0 B\big((k-r)^{N-j},(k-r+1)^{j-1}\big)_{\sh}\end{align*}
\proof
The proof follows by repeated applications of (\ref{un}), for the convenience of the reader we present the first 
step:
$$
\begin{array}{rcl}
B(k^N)&=&B(k^{N-1})_{\sh}\cup f_0 B\big((k-1)^1,k^{N-1}\big)\\
&=&B(k^{N-1})_{\sh}\cup f_0 B\big((k-1)^1,k^{n-2}\big)_{\sh}\cup f_0^2 B\big((k-1)^2,k^{N-2}\big)\\
&=&\ldots\\
&=&\bigcup_{r=0}^{N-1} f_0^r B\big((k-1)^{r},k^{N-1-r}\big)_{\sh}\cup f_0^NB\big((k-1)^N\big).
\end{array}
$$
The formula follows now by proceeding in the same way with $B\big((k-1)^N\big)$.
\endproof
\end{lem}

We will use Lemma~\ref{mainlem} to prove the following theorem
\begin{thm}\label{mainthm}\mbox{}
A $\PBW$ type basis of the truncated Weyl module $\W(kN,N)$ is given by 
$$
B(k^N)=\Big\{f^{i_0}_0f^{i_1}_1\cdots f^{i_{N-1}}_{N-1}\mid(i_0,\ldots,i_{N-1})\in S(k^N)\Big\}.
$$ 
\begin{example}\mbox{}
\begin{enumerate}
\item For $N=1$ we get that $S(k)$ is the set of $1$--tuples $(i_0)$ satisfying 
$$
i_0=\sum_{j=0}^0\frac{1!}{1-j}i_j\le 1! k - \sum_{\ell=0}^{-3}\frac{1!}{(1-\ell)!}(1-\ell-2)!b_\ell=k,
$$
so $S(k)=\{0,1,\ldots,k\}$ and $B(k)=\{f_0^{j}\mid j=0,\ldots,k\}$.

\item For $N=4$ and $k=2$ we get that $S(2^4)$ is the set of quadruples $(i_0,i_1,i_2,i_3)$ satisfying
$$
6 i_0+8 i_1+ 12 i_2 +24 i_3\le 48-2b_0,
$$
where $i_0=b_0 \bmod 4$ and 
$$
B(2^4)=\Big\{f_0^{i_0} f_1^{i_1}f_2^{i_2}f_3^{i_3}\mid (i_0,i_1,i_2,i_3)\in S(2^4)\Big\}.
$$
\end{enumerate}
\end{example}
\proof
The proof of Theorem~\ref{mainthm} proceeds by upward induction on $N$. The initial step is obvious (see also the previous example) and the induction begins.
So suppose that the theorem holds for all integers less than $N$. We claim the following:\vskip 4pt
\textit{Claim:}
For all $M<N$ we have
$$B\big(k^{M-j},(k+1)^{j}\big)=\Big\{f^{i_0}_0f^{i_1}_1\cdots f^{i_{M-1}}_{M-1}\mid(i_0,\ldots,i_{M-1})\in S\big(k^{M-j},(k+1)^{j}\big)\Big\}.$$

\vskip 3pt

\textit{Proof of the Claim:}
We prove this statement by induction. Note that there is nothing to prove if $j=0$. Hence we can suppose that $j>0$. We obtain 

\begin{align*}B\big(k^{M-j},(k+1)^{j}\big)=&B\big(k^{M-j},(k+1)^{j-1}\big)_{\sh}\cup f_0 B\big(k^{M-j+1},(k+1)^{j-1}\big)\\
=&\Big\{f^{i_0}_0f^{i_1}_1\cdots f^{i_{M-2}}_{M-2}\mid(i_0,\ldots,i_{M-2})\in S\big(k^{M-j},(k+1)^{j-1}\big)\Big\}_{\sh}
\\ \cup f_0&\Big\{f^{i_0}_0f^{i_1}_1\cdots f^{i_{M-1}}_{M-1} \mid (i_0,\ldots,i_{M-1})\in S\big(k^{M-j+1},(k+1)^{j-1}\big)\Big\}.\end{align*}

The shift by the map $\sh$ leads to the following description
\begin{align*}&\Big\{f^{i_0}_0f^{i_1}_1\cdots f^{i_{M-2}}_{M-2}\mid(i_0,\ldots,i_{M-2})\in S\big(k^{M-j},(k+1)^{j-1}\big)\Big\}_{\sh}
\\=&\Big\{f^{i_1}_1f^{i_2}_2\cdots f^{i_{M-1}}_{M-1}\mid\sum^{M-1}_{p=1}\frac{M!}{M-p}i_p\leq M!k-\sum^{M-4}_{\ell=1}\frac{M!}{(M-\ell)!}(M-\ell-2)!b_\ell+M(j-1)(M-2)!\Big\}\end{align*}
with
$$i_1-j+1=b_1\hspace{-0,3cm}\mod M-1,\ i_\ell+(b_{\ell-1}\hspace{-0,3cm}\mod M-\ell)=b_\ell\hspace{-0,3cm} \mod M-\ell,\ \mbox{ for } \ell=2,\dots,M-4$$
and 
\begin{align*}f_0&\Big\{f^{i_0}_0f^{i_1}_1\cdots f^{i_{M-1}}_{M-1}\mid(i_0,\ldots,i_{M-1})\in S\big(k^{M-j+1},(k+1)^{j-1}\big)\Big\}
\\=&\Big\{f^{i_0+1}_0f^{i_1}_1\cdots f^{i_{M-1}}_{M-1}\mid\sum^{M-1}_{p=0}\frac{M!}{M-p}i_p\leq M!k-\sum^{M-4}_{\ell=0}\frac{M!}{(M-\ell)!}(M-\ell-2)!b_\ell+(j-1)(M-1)!\Big\}
\\=&\Big\{f^{i_0}_0f^{i_1}_1\cdots f^{i_{M-1}}_{M-1}\mid\sum^{M-1}_{p=0}\frac{M!}{M-p}i_p\leq M!k-\sum^{M-4}_{\ell=0}\frac{M!}{(M-\ell)!}(M-\ell-2)!b_\ell+j(M-1)!,\ i_0\geq 1\Big\}\end{align*}
with
$$i_0-j=b_0\hspace{-0,3cm}\mod M,\ i_\ell+(b_{\ell-1}\hspace{-0,3cm}\mod M-\ell)=b_\ell\hspace{-0,3cm} \mod M-\ell,\ \mbox{ for } \ell=1,\dots,M-4.$$
Therefore, the claim follows with 
\begin{align*}&\Big\{f^{i_0}_0f^{i_1}_1\cdots f^{i_{M-2}}_{M-2}| (i_0,\ldots,i_{M-2})\in S\big(k^{M-j},(k+1)^{j-1}\big)\Big\}_{\sh}
\\=&\Big\{f^{0}_0f^{i_1}_1\cdots f^{i_{M-1}}_{M-1}\mid\sum^{M-1}_{p=0}\frac{M!}{M-p}i_p\leq M!k-\sum^{M-4}_{\ell=0}\frac{M!}{(M-\ell)!}(M-\ell-2)!b_\ell+j(M-1)!\Big\}.\end{align*}
Now it is easy to verify with Lemma~\ref{mainlem} that the theorem holds.
\endproof
\end{thm}
The proof of Theorem~\ref{mainthm} gives the following.
\begin{cor}\label{maincor}\mbox{}
A $\PBW$ type basis of the truncated Weyl module $\W(kN+j,N)$ is given by $$\Big\{f^{i_0}_0f^{i_1}_1\cdots f^{i_{N-1}}_{N-1}\mid(i_0,\ldots,i_{N-1})\in S\big(k^{N-j},(k+1)^{j}\big)\Big\}.$$
\end{cor}
\begin{rem}\label{rem1}
The fusion product $V(1)^{*N}$ is isomorphic to the truncated Weyl module $W_{\loc}(N,N)$ and also to the local Weyl module $W_{\loc}(N)$. The inductively obtained basis $B(1^N)$ coincides with the basis of the Weyl module $W_{\loc}(N)$ constructed in \cite{CP01}. However, we would like to emphasize that the PBW type basis of the truncated Weyl module $\W(m,N)$ described in Theorem~\ref{3} is different from the basis described in \cite[Section 6]{CV13}. For example, we have $f_1^3\in B(1^4)$ but $f_1^3$ is not contained in their basis.
\end{rem}

% \lambda^0\in\{\omega_1,\dots,\omega_n\} if $\lie g$ is of type $\tt A_n$
% \lambda^0\in\{\omega_1,\omega_n\} if $\lie g$ is of type $\tt B_n$
% \lambda^0\in\{\omega_1,\dots,\omega_n\} if $\lie g$ is of type $\tt C_n$
% \lambda^0\in\{\omega_1,\omega_{n-1},\omega_n\} if $\lie g$ is of type $\tt D_n$
% \lambda^0\in\{\omega_4\} if $\lie g$ is of type $\tt F_4$(\alpha_4 ist kurz)
% \lambda^0\in\{\omega_2\} if $\lie g$ is of type $\tt G_2$ (\alpha_2 ist kurz)
%%%%%%%%%%%%%%%%%%%%%%%%%%%%%%%%%%%%%%%%%%%%%
%         
%%%%%%%%%%%%%%%%%%%%%%%%%%%%%%%%%%%%%%%%%%%%%%%%%%%%%%%%%%%%%%%%%%%%%%%%%%%%%%%%%%%%%%%%%%%%%%%%%%%%%%%%%%%%%%%%%%%%%%%%%%%%%%%%%%%

% References
%%%%%%%%%%%%%%%%%%%%%%%%%%%%%%%%%%%%%%%%%%%%%%%%%%%%%%%%%%%%%%%%%%%
\bibliographystyle{plain}
\bibliography{fusionandtoroidal-biblist}
\end{document}